\documentclass[11pt]{article}
\usepackage[final]{epsfig}
\usepackage{graphics}
\usepackage{amsmath}
\usepackage{amsfonts}
\usepackage{latexsym}
\usepackage{amssymb}
\usepackage{graphicx}
\newtheorem{lemma}{Lemma}[section]
\newtheorem{proposition}[lemma]{Proposition}
\newtheorem{remark}[lemma]{Remark}
\newtheorem{example}[lemma]{Example}
\newtheorem{theorem}[lemma]{Theorem}
\newtheorem{definition}{Definition}

\begin{document}
\newcommand{\eps}{{\varepsilon}}
\newcommand{\proofend}{\hfill$\Box$\bigskip}
\newcommand{\C}{{\mathbb C}}
\newcommand{\Q}{{\mathbb Q}}
\newcommand{\R}{{\mathbb R}}
\newcommand{\T}{{\mathbb T}}
\newcommand{\Z}{{\mathbb Z}}
\newcommand{\RP}{{\mathbb {RP}}}
\newcommand{\A}{{\mathbb {A}}}

\def\proof{\paragraph{Proof.}}

%%%%%%%%%%%%%%%%%%%%%%%%%%%%%%%%%%%%%%%%%%%%%%
%%%%%%%%%%%%%%%%%%%%%%%%%%%%%%%%%%%%%%%%%%%%%%

\newcommand{\marginnote}[1]
{%\mbox{}\marginpar{\center{\hspace{0pt}\tiny{\bf#1}}}
}

\newcounter{bk}
\newcommand{\bk}[1]
{\stepcounter{bk}$^{\bf BK\thebk}$%
\footnotetext{\hspace{-3.7mm}$^{\blacksquare\!\blacksquare}$
{\bf BK\thebk:~}#1}}

\newcounter{st}
\newcommand{\st}[1]
{\stepcounter{st}$^{\bf ST\thest}$%
\footnotetext{\hspace{-3.7mm}$^{\blacksquare\!\blacksquare}$
{\bf ST\thest:~}#1}}

%%%%%%%%%%%%%%%%%%%%%%%%%%%%%%%%%%%%%%%%%%%%%%%%%%%%%%%%%%

\title {Contact complete integrability}
\author{Boris Khesin\thanks{
Department of Mathematics,
University of Toronto, Toronto, ON M5S 2E4, Canada;
e-mail: \tt{khesin@math.toronto.edu}
}
\, and Serge Tabachnikov\thanks{
Department of Mathematics,
Pennsylvania State University, University Park, PA 16802, USA;
e-mail: \tt{tabachni@math.psu.edu}
}
\\
}
\date{October 1, 2009}
\maketitle
\begin{abstract}
Complete integrability in a symplectic setting means the existence of a Lagrangian foliation leaf-wise preserved by the dynamics. In the
paper we describe complete integrability in a contact set-up as  a more subtle structure:  a flag of two foliations, Legendrian and co-Legendrian,
and a holonomy-invariant transverse measure of the former in the latter.  This  turns out to be equivalent to the existence of a
canonical $\R\ltimes \R^{n-1}$ structure on the leaves of the co-Legendrian foliation.   Further, the above structure implies the existence of $n$ contact fields preserving a special contact 1-form, thus providing the geometric framework and  establishing equivalence with previously known definitions of contact integrability. 
We also show that contact completely integrable systems are solvable in quadratures.

We present an example of contact complete integrability: the billiard system inside an ellipsoid in pseudo-Euclidean space,  restricted to the space of oriented null geodesics. We describe a surprising acceleration mechanism for closed light-like billiard trajectories. 
\end{abstract}

\tableofcontents

\section{Introduction} \label{intro}

Our first motivation is the following V. Arnold's problem No 1995--12 in \cite{Arnold}:
\begin{quote}
Transfer the theory of completely integrable Hamiltonian systems from symplectic geometry to contact geometry (where, e.g., the Lagrangian invariant manifolds with their natural affine structures determined by Lagrangian fibrations must be substituted by Legendrian invariant manifolds with their natural projective structures determined by Legendrian fibrations). Carry over the Liouville theorem to this context and find applications to the infinite-dimensional case (where the equations of characteristics are partial differential).
\end{quote}

The classical set-up for the Arnold-Liouville theorem is a symplectic manifold $(M^{2n},\omega)$  (for example, the phase space of a mechanical system) and a discrete- or continuous-time symplectic dynamical system on it, that is, a symplectomorphism $T:M\to M$  or a symplectic vector field $v$ on $M$, respectively. 
(Here and elsewhere we refer to \cite{A-G} for a succinct exposition of the basic facts of symplectic and contact geometry.) 
Recall that a Lagrangian manifold $F^n\subset M^{2n}$ is a half-dimensional submanifold such that the restriction of $\omega$ to $F$ vanishes.  
A symplectic dynamical system is called completely integrable if $M$ is endowed with a Lagrangian foliation ${\cal F}$ whose leaves are invariant under the dynamics.

A fundamental geometrical fact, underlying the Arnold-Liouville theorem, is that the leaves of a Lagrangian foliation carry a canonical affine structure. Choose $n$ functionally independent ``integrals" (functions, constant on the leaves of  ${\cal F}$) and consider their symplectic gradients. One obtains $n$ commuting vector fields, tangent to the leaves of  ${\cal F}$ and providing a field of frames along the leaves. A different choice of integrals results in applying  linear transformations, constant along each leaf, to these frames. Thus each leaf has a flat structure.

 The map $T$, or the vector field $v$, preserves the symplectic structure and the foliation leaf-wise, and hence preserves the affine structure on the leaves. It follows that $T$ is a parallel translation, and $v$ is a constant vector field, on each leaf of the Lagrangian foliation. 
 
This has strong dynamical consequences. If a point  is periodic then so are all the points on the same  leaf of ${\cal F}$, and with the same period. This implies Poncelet-style theorems (see, e.g., \cite{L-T} for a recent application). If a leaf is compact, it must be a torus, and the dynamics is a quasi-periodic motion on the torus. Another useful consequence: if two symplectic maps share an invariant Lagrangian foliation then they commute (because so do parallel translations). We refer to \cite{Ve1} for more detail on complete integrability in the discrete-time case.

To summarize, the definition of a completely integrable dynamical system consists of two parts: a certain geometrical structure on a symplectic manifold $M$, namely, a Lagrangian foliation ${\cal F}$, and a discrete- or continuous-time symplectic dynamical system, preserving this structure. It is natural  to call the first part, the pair $(M,{\cal F})$, a completely integrable symplectic manifold.

Contact manifolds are odd-dimensional relatives of symplectic manifolds. Let $(M^{2n-1},\xi)$ be a contact manifold with a contact distribution $\xi$. Recall that a Legendrian submanifold  $F^{n-1}\subset M^{2n-1}$ is an integral manifold of $\xi$ of the maximal possible dimension $n-1$. The leaves of a Legendrian foliation carry a canonical projective structure: this is a contact counterpart to the above-mentioned theorem about Lagrangian foliations (we shall dwell on this projective structure in Sections \ref{interp}). The problem is to extend the notion of complete integrability to contact manifolds.

Note that the simplest particular case of a contact manifold is a 1-dimensional manifold, $\R^1$ or $S^1$, with the trivial contact structure. A natural definition of  integrability in dimension 1 (with discrete- or continuous-time) is the existence of a non-vanishing invariant differential 1-form; for a diffeomorphism of $S^1$, this implies that the map is conjugated to a rotation.

Contact complete integrability was studied before: see \cite{B-M1,B-M2,Ba,Ler,Mir} and also related papers on Legendrian foliations \cite{Pang,Lib}. For example, according to \cite{B-M1}, a completely integrable vector field on a  contact manifold $M^{2n-1}$ is the Reeb field of a contact 1-form, for which the space of
first integrals determines a fibration with $n$-dimensional fibers defined locally by
the action of a torus $\T^n$ of contact transformations, 
%\st{This is taken from the MathSciNet review of \cite{B-M1}} 
see Section \ref{earlier} for a brief survey of earlier work.

The main goal in this paper is to give a definition that extends the  earlier ones and that works equally well in the continuous- and discrete-time cases. 
Our second motivation was to place the recently studied examples \cite{Gen-Kh-Tab, Kh-Tab} into the general context of contact complete integrability.

These examples include the geodesic flow 
on an ellipsoid and the billiard map inside an ellipsoid in pseudo-Euclidean space. In pseudo-Euclidean setting, one has a trichotomy for an oriented line: it may be space-like, time-like, and light-like (or null), that is,  having positive, negative, or zero energy. The manifolds of oriented non-parameterized space- and time-like lines carry canonical symplectic structures, just like in the Euclidean case, but the space of null lines has a canonical contact structure; see \cite{Kh-Tab} and  Section \ref{null}.

Let $S$ be a smooth closed hypersurface in a pseudo-Euclidean space. The billiard system inside $S$ can be considered as a map on the space of oriented lines taking the incoming billiard trajectory to the outgoing one. The law of reflection is determined by the energy and momentum conservation, therefore the type of a line (space-, time-, or light-like) does not change. Restricted to space- and time-like lines, the billiard transformation is a symplectic map, but its restriction to light-like lines is a contact map. 

%Likewise, the geodesic flow on $S$ can be viewed as a vector field on the space of oriented lines tangent to $S$; this flow preserves the type of a line too. Restricted to space- and time-like lines, the geodesic flow is a Hamiltonian vector field, but its restriction to light-like lines is a contact vector field. \st{Is this so? Razobrat'sya!}

If $S$ is an ellipsoid, the respective billiard transformation is integrable, in the following sense. An ellipsoid in $n+1$-dimensional pseudo-Euclidean space determines a pseudo-confocal family of quadrics, see \cite{Kh-Tab} and  Section \ref{null}. A space- or time-like billiard trajectory remains tangent to $n$ fixed pseudo-confocal quadrics. This gives $n$ integrals of the billiard map on the $2n$-dimensional symplectic spaces of oriented space- or time-like lines. These integrals Poisson commute and hence define an invariant  Lagrangian foliation. This is just like  the Euclidean case, see, e.g., \cite{Ta1,Ta2}.

However, we lose one integral on the space of null lines: a light-like billiard trajectory remains tangent to $n-1$ fixed pseudo-confocal quadrics. This gives $n-1$ integrals on the $2n-1$-dimensional contact space of oriented light-like lines and hence a foliation ${\cal F}^{n}$. It turns out that the distribution given by the intersection of the tangent spaces to the leaves of ${\cal F}$ with the contact hyperplanes is also integrable, and one obtains a Legendrian foliation ${\cal G}^{n-1}$ whose leaves foliate the leaves of ${\cal F}$. Furthermore,  the billiard transformation has an invariant contact form -- morally, another integral, since all contact forms for a given contact structure differ by multiplication by a non-vanishing function -- and this additional integral commutes, in an appropriate sense, with the other $n-1$ integrals, see Section \ref{null}. 

\medskip

The above motivates the following general definition of contact integrability. 
Let $(M^{2n-1},\xi)$ be a contact manifold with contact distribution $\xi$. A foliation ${\cal F}^{n}$ is called {\it co-Legendrian} if it is transverse to $\xi$ and the distribution $T{\cal F} \cap \xi$ is integrable. Let  ${\cal G}^{n-1}$ be the respective Legendrian foliation. We have a flag of foliations $({\cal F},{\cal G})$. In Section \ref{coLeg} we show that the canonical projective structure on the leaves of ${\cal G}$ reduces to an affine structure.

If a contact dynamical system preserves a co-Legendrian foliation leaf-wise then it sends the leaves of the respective Legendrian foliation ${\cal G}$ to each other, preserving the affine structures therein. Thus the dynamics reduces to 1-dimensional one on the space of leaves of ${\cal G}$ within a leaf of ${\cal F}$. For this dynamics to be integrable, one needs an invariant 1-form on this 1-dimensional space of leaves. 

\begin{definition} \label{main}
{\rm A {\it completely integrable contact manifold} 
$M$ is a co-Legendrian foliation ${\cal F}$ on $M$ such that, for each leaf $F$ of ${\cal F}$, the respective codimension one foliation ${\cal G}$ on $F$ has a holonomy invariant transverse smooth measure. A discrete- or continuous-time {\rm contact completely integrable system} 
on $M$ is a contactomorphism, or a contact vector field, that preserves ${\cal F}$  leaf-wise and preserves the above transverse measure of the foliation ${\cal G}$.}
\end{definition} 

We show in Section \ref{semsect} that the leaves of a co-Legendrian foliation on a completely integrable contact manifold have a canonical $\R\ltimes\R^{n-1}$-structure. This has strong dynamical implications, similarly to the  flat $\R^n$-structure on the leaves of a Lagrangian foliation of a symplectic manifold.

An example of a completely integrable contact manifold $M$ is analyzed in  Sections \ref{contf} 
and \ref{commf}: $M$ has a contact form whose Reeb field is tangent to the co-Legendrian foliation. We show that then the contact form defines a holonomy invariant transverse smooth measure on the space of leaves of the Legendrian foliation within a leaf of the co-Legendrian one. We show in Section \ref{intgeo} that the familiar example of a completely integrable geodesic flow on a Riemannian manifold fits into this framework.

In Section \ref{null}, we show that the billiard ball map inside an ellipsoid in pseudo-Euclidean space, restricted to oriented light-like lines, is a completely integrable contact map. We do this by constructing an invariant contact form on the contact space of oriented null lines whose Reeb field is tangent to a co-Legendrian foliation.

%%%%%%%%%%%%%%%%%%%%%%%%%%%%%%%%%%%%%%%%%%%%%%%%%%%%%
%%%%%%%%%%%%%%%%%%%%%%%%%%%%%%%%%%%%%%%%%%%%%%%%%%%%%

\section{Geometry of co-Legendrian foliations} \label{coLeg}

In this section we study the geometry of co-Legendrian foliations. 

\subsection{Example of a co-Legendrian foliation} \label{examples} 

\begin{example} \label{ex}
{\rm Let $M$ be a contact manifold with a Legendrian foliation ${\cal G}$, and let $\phi_t$ be a 1-parameter group of contactomorphisms preserving this foliation. Assume that the vector field corresponding to $\phi_t$ is transverse to the contact distribution. Then, acting by $\phi_t$ on ${\cal G}$, yields a co-Legendrian foliation ${\cal F}$, that is, the leaves of ${\cal F}$ are the orbits of the leaves of ${\cal G}$ under the flow $\phi_t$.}
\end{example}

In fact, this example is universal, as the following lemma shows.

\begin{lemma} \label{canon}
Every co-Legendrian foliation is locally contactomorphic to the one in Example \ref{ex}.
\end{lemma}

\proof Recall that a contact element on a smooth manifold $N$ is a hyperplane in a tangent space to $N$. Since a contact element is the kernel of a covector, uniquely determined up to a non-zero factor, the space of contact elements is $PT^*N$, the projectivization of the cotangent bundle. The space of contact elements has a canonical contact structure given by the ``skating condition": the velocity of the foot point of a contact element lies in this contact element. One has the fibration $p: PT^*N\to N$ whose fibers are Legendrian submanifolds; these fibers consist of contact elements with a fixed foot point. Every Legendrian foliation is locally contactomorphic to this one.

Suppose that $M$ is a contact manifold, ${\cal F}$ is a co-Legendrian foliation and ${\cal G}$ the respective Legendrian foliation. We may assume that, locally, $M=PT^*N$ and ${\cal G}$ is the fibration $p: PT^*N\to N$. Then  ${\cal F}$ projects to  a 1-dimensional foliation ${\cal L}$  in $N$. In other words, a leaf of ${\cal F}$ consists of contact elements whose foot points lie on a leaf of ${\cal L}$.

Consider a 1-parameter group of diffeomorphisms of $N$ whose trajectories are the leaves of ${\cal L}$. Diffeomorphisms of $N$ naturally act on contact elements on $N$, so we obtain a 1-parameter group of contactomorphisms of $M$ preserving the foliation ${\cal G}$. If we restrict to the open set of contact elements on $N$ that are not tangent to ${\cal L}$ then ${\cal F}$ is obtained from ${\cal G}$ as in Example \ref{ex}.
\proofend

%%%%%%%%%%%%%%%%%%%%%%%%%%%%%%%%%%%%%%%%%%%%%%%%%%%%%

\subsection{Symplectic interpretation of co-Legendrian foliations} \label{interp}

Let $(M^{2n-1},\xi)$ be a contact manifold. Recall the notion of symplectization (or the symplectic cone). Let $P^{2n}\subset T^*M$ be the total space of the principle $\R^*$-bundle $\pi :P\to M$ whose fibers consist of non-zero covectors $(x,p)$ that vanish on the contact element $\xi(x)$ in $M$ at their respective foot points $x$. The symplectization $P$ has a canonical 1-form $\lambda$, the restriction of the Liouville 1-form in $T^*M$, and the 2-form $\omega=d\lambda$ is a symplectic structure on $P$. One has the multiplicative $\R^*$-action on $P$; let $E$ be the respective  vector field, called the Euler field. The following identities hold:
\begin{equation} \label{sy}
i_E\omega=\lambda,\ \ \lambda(E)=0,\ \  L_E \lambda=\lambda.
\end{equation}

For example, the symplectization of the space of contact elements $PT^*N$ is the cotangent bundle $T^*N$ with the zero section removed. The $\R^*$-action is the fiber-wise scaling of covectors, and the 1-form $\lambda$ is the Liouville form in $T^*N$.

The operation of symplectization relates the contact geometry of $M$ to the homogeneous symplectic geometry of  $P$. Specifically, contactomorphisms of $M$ are the symplectomorphisms of $P$ that commute with the $\R^*$-action; the preimage of a Legendrian submanifold in $M$ is an $\R^*$-invariant (conical) Lagrangian submanifold in $P$;  the preimage of a Legendrian foliation in $M$ is an $\R^*$-invariant Lagrangian foliation in $P$, etc.

Let ${\cal F}$ be a co-Legendrian foliation on $M$ and ${\cal G}$ the respective Legendrian foliation. Set: $\overline {\cal F}=\pi^{-1} ({\cal F}),\ \overline {\cal G}=\pi^{-1} ({\cal G})$. In the next lemma, we interpret co-Legendrian foliations in symplectic terms. 

\begin{figure}[hbtp]
\centering
\includegraphics[width=3.5in]{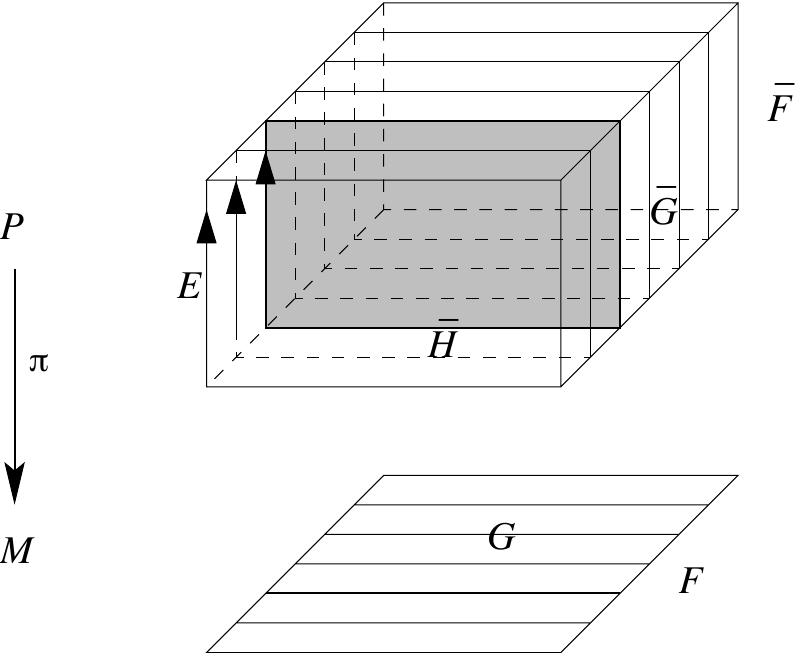}
\caption{A leaf $F$ of the foliation $\overline {\cal F}$ in $P$, foliated by $\overline {\cal G}$ and by $\overline {\cal H}$  projects to a leaf $F$ of the foliation ${\cal F}$ in M, foliated by  ${\cal G}$}
\label{leaves}
\end{figure}

\begin{lemma}
$\overline {\cal F}^{n+1}$ is a co-isotropic foliation in $P$. Its symplectic orthogonal complement $\overline {\cal H}^{n-1}$ is an isotropic foliation transverse to $E$, and  $\overline {\cal G}^n$ is  spanned by $E$ and $\overline {\cal H}^{n-1}$ (see figure \ref{leaves}). Conversely, given  a co-isotropic foliation $\overline {\cal F}^{n+1}$ in $P$, tangent to the Euler field $E$ and transverse to $\ker \lambda$, the projection of $\overline {\cal F}$ to $M$ is a  co-Legendrian foliation therein. 
\end{lemma}

\proof Let $f_i:M\to \R,\ i=1,\dots,n-1$, be locally defined functions whose common level surfaces are the leaves of ${\cal F}$, and let $\bar f_i=\pi^*(f_i)$. Then the homogeneous functions $\bar f_i:P\to \R$ of degree zero define the foliation $\overline {\cal F}$. 

First, we show that the symplectic orthogonal complement to the tangent space $T\overline {\cal F}$ is spanned by the Hamiltonian vector fields sgrad $\bar f_i$.  Indeed, 
consider a vector $v\in T\overline {\cal F}$. Then $\omega({\rm sgrad}\bar f_i, v)=d \bar f_i (v) =0,$
since $\bar f_i$ is constant on the leaves of $\overline {\cal F}$. 

Next, we show that the distribution spanned by the Hamiltonian vector fields sgrad $\bar f_i$ is integrable. Indeed, this distribution is isotropic, hence 
$\omega({\rm sgrad}\bar f_i, {\rm sgrad}\bar f_j)=0=\{f_i, f_j\}$. It follows that $[{\rm sgrad}\bar f_i, {\rm sgrad}\bar f_j]=0$, so $\overline {\cal H}$ is a foliation.

One has: $\overline {\cal G}=\overline {\cal F} \cap \ker \lambda$. We claim that $\overline {\cal G}$ is spanned by $E$ and ${\rm sgrad} \bar f_i$. Indeed, $E$ and ${\rm sgrad} \bar f_i$ are tangent to $\overline {\cal F}$. One has: $\lambda(E)=0$ and 
$$\lambda ({\rm sgrad} \bar f_i)=\omega(E,{\rm sgrad}\bar f_i)=-d\bar f_i(E)=0,$$
since $\bar f_i$ is homogeneous of degree zero with respect to the Euler field. Thus $E$ and all ${\rm sgrad} \bar f_i$ are tangent to $\overline {\cal G}$.

Let us check that $E$ is transverse to  $\overline {\cal H}$. If not, then, at some point,  $E={\rm sgrad} \bar f$ for a function $f:M\to \R$ that is constant on the leaves of ${\cal F}$. Then at that point
$
\lambda=i_E\omega=i_{{\rm sgrad} \bar f}\omega=d \bar f.
$
This is a contradiction since the foliation ${\cal F}$ is transverse to the contact structure, and hence $\lambda$ does not vanish on the tangent spaces to its leaves.

Finally, we claim that if $\overline {\cal F}^{n+1}$ is a co-isotropic foliation in $P$, tangent to the Euler field $E$, then the projection $\pi :P\to M$ takes $\overline {\cal F}$ to a co-Legendrian foliation. Indeed, the foliation $\overline {\cal F}$ is invariant under the Euler field since $E$ is tangent to it. Thus  $\overline {\cal F}$ is conical. Then the distribution $T\overline {\cal F} \cap \ker \lambda$ is a conical Lagrangian foliation that projects to a Legendrian foliation in $M$.
\proofend

Thus a co-Legendrian foliation on a contact manifold $M^{2n-1}$ is the same as a co-isotropic $n+1$-dimensional foliation on its symplectization $P^{2n}$ given by $n-1$ Poisson commuting homogeneous functions of degree zero.

%%%%%%%%%%%%%%%%%%%%%%%%%%%

\subsection{Flat structure on the leaves of ${\cal G}$} \label{affsect}

As we mentioned in Introduction, the leaves of a Legendrian foliation carry a canonical projective structure. Let us recall this construction. 

Let $(M,\xi)$ be a contact manifold and ${\cal G}$ a Legendrian foliation. As before, we may assume that  $M=PT^*N$ and ${\cal G}$ is the fibration $p: PT^*N\to N$. 
Let $x\in N$ and $G_x=p^{-1}(x)$. Then $dp$ takes the contact hyperplanes along the leaf $G_x$ to hyperplanes in the tangent space $V:=T_x N$. The set of all such
hyperplanes  is $P(V)=\RP^{n-1}$, and we obtain a mapping $\varphi: G_x\to \RP^{n-1}$. Due to complete non-integrability of the contact structure, $\varphi$ is a local diffeomorphism. Thus $G_x$ has a projective structure.

Now let $(M,\xi)$ be a contact manifold, and ${\cal F}$ and ${\cal G}$ be co-Legendrian and the respective Legendrian foliations.

\begin{lemma} \label{affstr}
The projective structure on the leaves of ${\cal G}$ has a reduction  to an affine structure.
\end{lemma}

\proof  In the notation of the preceding paragraphs,
the tangent spaces to a leaf of ${\cal F}$ are taken by $dp$ to a line $\ell\subset V$.  The set of hyperplanes in $V$ passing through $\ell$ is a projective hyperplane 
$\RP^{n-2} \subset P(V)=\RP^{n-1}$, and the image of $\varphi$ does not intersect this 
projective hyperplane. The complement $\RP^{n-1} - \RP^{n-2}$ is an affine space.
Thus we have a local diffeomorphism $\varphi: G_x\to \A^{n-1}$, whence an affine structure on the leaves of ${\cal G}$.
\proofend

As usual, the existence of an affine structure imposes restrictions on the topology of the leaves. For example, a compact leaf of ${\cal G}$ is a torus.

\begin{remark} \label{aff}
{\rm Alternatively, one can define an affine structure on the leaves of ${\cal G}$ as follows. Recall that the isotropic foliation $\overline {\cal H}$ is generated by the vector fields ${\rm sgrad}\bar f_i$, where the functions $\bar f_i$ are homogeneous of degree zero. The commuting vector fields ${\rm sgrad}\bar f_i$ define an affine structure on the leaves of $\overline {\cal H}$.

Since the functions $\bar f_i$ are homogeneous of degree zero,  $[E,{\rm sgrad}\bar f_i]=-{\rm sgrad}\bar f_i$. Therefore the $\R^*$-action preserves the foliation $\overline {\cal H}$, sending leaves to leaves, and these maps preserve the affine structure on the leaves. The projection $\pi: P\to M$ diffeomorphically maps the leaves of $\overline {\cal H}$ to the leaves of the Legendrian foliation ${\cal G}$ endowing the latter with an affine structure. 
}
\end{remark}

%%%%%%%%%%%%%%%%%%%%%%%%%%%%%%%%%%%%%%%%%%%%%%%%%%%%%%%%%%%%%%%%%%%%%%%%%%

\subsection{Weakly integrable contact systems} \label{weaksect}

\begin{definition} \label{weak}
{\rm A discrete- or continuos-time contact {\it weakly integrable system} is a contact dynamical system on a contact manifold that has a leaf-wise invariant co-Legendrian foliation.}
\end{definition}

Such a system reduces to a 1-dimensional one. The leaves of the Legendrian foliation ${\cal G}$ within a leaf $F$ of the co-Legendrian foliation ${\cal F}$ are mapped to each other by parallel translations in their respective affine coordinates, but the motion on the 1-dimensional space of leaves ${F}/{\cal G}$ may be arbitrary. 
As Lemma \ref{canon} shows, any diffeomorphism of $N$, preserving the one-dimensional foliation ${\cal L}$ leaf-wise, lifts to a weakly integrable contactomorphism of (an open subspace of) the space of contact elements of $N$.

If $M$ is 1-dimensional, the co-Legendrian foliation consists of one leaf, $M$ itself, and 
the definition imposes no constraints on the dynamics.

%%%%%%%%%%%%%%%%%%%%%%%%%%%%%%%%%%%%%%%%%%%%%%%%%%%%%%%%%%%%%%%%%%%%%%%%%%
%%%%%%%%%%%%%%%%%%%%%%%%%%%%%%%%%%%%%%%%%%%%%%%%%%%%%%%%%%%%%%%%%%%%%%%%%%

\section{Completely integrable contact manifolds} \label{other}

In this section we study the geometry of completely integrable contact manifolds and completely integrable contact dynamical systems.

\subsection{Semi-direct product structure} \label{semsect}

Let $G$ be a subgroup of the group of diffeomorphisms of $\R^n$. A $G$-structure on an $n$-dimensional manifold is (an equivalent class of) an atlas whose transition functions belong to $G$. In these terms, the leaves of a Lagrangian foliation have an $\R^n$-structure, where $\R^n$ is the group of parallel translations of $n$-dimensional affine space.

Let $M^{2n-1}$ be a completely integrable contact manifold with the flag of co-Legendrian and Legendrian foliations $({\cal F}^n, {\cal G}^{n-1})$. Let $\R\ltimes\R^{n-1}$ be a semi-direct product of $\R$ and $\R^{n-1}$:
$$
0\to \R^{n-1} \to \R \ltimes \R^{n-1} \to \R \to 0.
$$

\begin{example} 
{\rm Given a number $\lambda\in \R$ and a vector $b\in \R^{n-1}$ consider  affine maps $v\mapsto e^\lambda v +b$ of the space $\R^{n-1}\ni v$. Then  the set of such pairs $(\lambda,b)$ forms a Lie group with respect to natural composition of the affine maps. 
Similarly, one can define the Lie group by composing affine transformations $v\mapsto e^\lambda Pv +b$ 
for a projector $P: \R^{n-1}\to\R^{n-1}$ (with $P^2=P$). These Lie groups give examples of  semi-direct products $\R\ltimes \R^{n-1}$. The first example corresponds to the case of $P=id$, while the 
direct product group $\R^n=\R\times\R^{n-1}$ corresponds to $P=0$.}
\end{example}

\begin{lemma} \label{sem}
The leaves of the co-Legendrian foliation ${\cal F}$ of a completely-integrable contact manifold have a canonical $\R\ltimes\R^{n-1}$-structure.
\end{lemma}

\proof
Let $F$ be a leaf of ${\cal F}$. One has an exact sequence of vector bundles:
$$
0\to T{\cal G} \to TF \to N{\cal G}=TF/T{\cal G} \to 0
$$
($N{\cal G}$ is the normal bundle of the foliation ${\cal G}$).
By Lemma \ref{affstr}, the leaves of ${\cal G}$ have an $\R^{n-1}$-structure. The transverse invariant measure of the foliation ${\cal G}$ in $F$ fixes a trivialization of $N{\cal G}$. The two combined yield a $\R\ltimes\R^{n-1}$-structure on $F$. 
\proofend

Recall that a completely integrable system on a symplectic manifold $M^{2n}$ can be defined by a local $\R^n$-action which preserves the symplectic structure and is generically free.

Similarly, a completely integrable system on a contact manifold can be defined by a $\R\ltimes\R^{n-1}$-action, where the abelian subgroup $\R^{n-1}$ acts locally free along the contact planes. 
Namely, consider a contact manifold $M^{2n-1}$ with a contact distribution $\xi$. 
One can see that Lemma \ref{sem} is equivalent to the following

\begin{lemma} \label{sem-action}
The existence of a co-Legendrian foliation with an invariant measure on a contact manifold $(M, \xi)$ is equivalent to the existence of a local  $\R\ltimes\R^{n-1}$-action on $M$ such that the $\R^{n-1}$-orbits of the abelian subgroup are tangent to the distribution $\xi$.
\end{lemma}

By construction, the $\R\ltimes\R^{n-1}$-orbits define leaves of the co-Legendrian foliation ${\cal F}$, while the orbits of the abelian subgroup provide the  Legendrian foliation ${\cal G}$. The $\R$-action in the semi-direct product gives the holonomy-invariant transverse measure. Conversely, the existence of the $\R\ltimes\R^{n-1}$-structure for
${\cal F}$ implies the existence of a local $\R\ltimes\R^{n-1}$-action in each leaf of 
${\cal F}$, whose $\R^{n-1}$-orbits are leaves of the Legendrian foliation ${\cal G}$.

%%%%%%%%%%%%%%%%%%%%%%%%%%%%%%%%%%%%%%%%%%%%%%%%%%%%%%%%%%%%%%

\subsection{A special contact form} \label{contf}

It turns out that the existence of a local $\R\ltimes\R^{n-1}$-action on a contact manifold $M$ implies the
existence of a special contact form whose Reeb field is tangent to the group orbits.
Recall that the {\it Reeb vector field} $v$ of a contact form $\lambda$ spans the kernel of $d\lambda$ and is normalized by the condition $\lambda(v)=1$.

\begin{lemma} \label{inv-form}
There is a contact form $\lambda$ on $M$ whose Reeb field is tangent to the $\R\ltimes\R^{n-1}$-orbits.
\end{lemma} 

\proof 
Let $V\in (\R\ltimes\R^{n-1})$ be a generic element of the Lie algebra and $U\in \R^{n-1}$ an element of the abelian subalgebra. Let $v$ and $u$ be the corresponding vector fields on $M$.
Note that in the semi-direct product $\R\ltimes\R^{n-1}$ the elements $U$ and $V$ (and hence the fields $u$ and $v$) satisfy the relation $[v,u]=au$ for some $a\in \R$. Also note that by the definition of the action, $u$ is tangent to the distribution $\xi$, while $v$ is transversal to $\xi$ at a generic point. 

Now let $\lambda_0$ be a contact 1-form defining the contact structure $\xi$. Define the 1-form $\lambda$
by normalizing $\lambda_0$ as follows:  $\lambda=\lambda_0/\lambda_0(v)$, so that $\lambda(v)=1$.\footnote{This construction of the invariant 1-form is similar to the one for the action of an abelian group discussed in \cite{B-M2}.}

We see that, on the one hand, 
$$
i_{[v,u]}\lambda=i_{au}\lambda=0\,.
$$

On the other hand, 
$$
i_{[v,u]}\lambda=i_v L_u\lambda-L_u i_v\lambda= i_v i_u d\lambda+ i_v d i_u \lambda  -L_u i_v\lambda=i_v i_u d\lambda\,,
$$
where we used that $i_v\lambda=1$ and $i_u\lambda=0$.
The equality $ i_v i_u d\lambda=0$ implies that the kernel of the 2-form $d\lambda$ is tangent to the orbits of the $\R\ltimes\R^{n-1}$-action. Indeed, this equality shows that the projection of $v$ to planes of $\xi$ along this kernel must be tangent to the Legendrian $\R^{n-1}$-orbits. Thus the Reeb field for the 1-form $\lambda$ belongs to the $\R\ltimes\R^{n-1}$-orbits.
\proofend

The following lemma shows that the converse statement also holds: the existence of such a special contact 1-form is equivalent to the existence of co-Legendrian foliation with a transverse measure. 

Let $(M,\xi)$ be a contact manifold with a contact form $\lambda$, and let ${\cal F}$ and ${\cal G}$ be a co-Legendrian and the respective Legendrian foliations. Does a contact form $\lambda$ determine  a holonomy invariant transverse smooth measure of the foliation ${\cal G}$ within the leaves of ${\cal F}$? The next lemma also gives a sufficient condition.

\begin{lemma} \label{Reeb}
Assume that the Reeb field of $\lambda$ is tangent to the co-Legendrian foliation  ${\cal F}$. Let $F$ be a leaf of  ${\cal F}$. Then the 1-form $\lambda$ determines a holonomy invariant transverse smooth measure of the foliation ${\cal G}$ on the manifold $F$.
\end{lemma}

\proof
We need to check that the restriction of $\lambda$ to $F$ is a basic differential form with respect to the foliation ${\cal G}$; this means that for every vector field $u$, tangent to ${\cal G}$, one has: $i_u\lambda=L_u\lambda=0$. If $\lambda$ is basic then it descends on the (locally defined) space of leaves and defines a 1-form on this space.

We have $\lambda(u)=0$, since ${\cal G}$ is Legendrian. Then, by Cartan's formula,  $L_u\lambda=i_u d\lambda$, and we want to show that $i_u d\lambda=0$. The tangent space $TF$ is spanned by $T{\cal G}$ and $v$, the Reeb field. If $w\in T{\cal G}$ then $d\lambda(u,w)=0$ since $T{\cal G}$ is a Lagrangian subspace of the symplectic space  $\xi=\ker \lambda$ with the symplectic structure $d\lambda$. On the other hand, $d\lambda(u,v)=0$ since $v \in \ker \lambda$. Thus $i_u d\lambda=0$, and we are done.
\proofend

The above two lemmas give a necessary and sufficient condition of contact integrability in terms of a special 1-form.

\begin{theorem} \label{nec-suff}
The existence of a co-Legendrian foliation  ${\cal F}$ with an invariant transverse smooth measure is equivalent to the existence of a co-Legendrian foliation with a special contact 1-form whose Reeb field is tangent to the foliation.
\end{theorem}

%%%%%%%%%%%%%%%%%%%%%%%%%%%%%%%%%%%%%%%%%%%%%%%%%%%%%%%%%%%%%%

\subsection{Commuting fields and invariant contact forms} \label{commf}

The existence of a local $\R\ltimes\R^{n-1}$-action also implies the existence of an appropriate local $\R^n$-action (with the same orbits), preserving this form: one can define $n$ commuting vector fields which leave the contact form $\lambda$ invariant and span the same foliation ${\cal F}$.
Note however, that although the $\R\ltimes\R^{n-1}$- and  $\R^n$-orbits coincide, the orbits of the $\R^{n-1}$-subgroups in these two groups are different: any nonzero vector field preserving a contact form cannot be tangent to a contact distribution. 

Recall that if a contact form $\lambda$ is chosen on a contact manifold $(M,\xi)$ then one can assign a contact vector field $X_f$ to a smooth function $f$: a contact form determines a section of the symplectization $P$,   and this makes it possible to extend  $f$ to $P$ as a homogeneous of degree one function; the Hamiltonian vector field of this extended function projects to a contact vector field $X_f$ on $M$. The correspondence between the functions and contact vector fields is described by the formula $\lambda(X_f)=f$. In particular, for $f\equiv 1$ one has $X_1=v$, the Reeb field. Note also the formula: $L_{X_f}\lambda=df(v) \lambda$.

Further, one defines the Jacobi bracket on smooth functions: $[f,g]=\lambda([X_f,X_g])$. This operation satisfies the Jacobi identity, but not the Leibniz one. The correspondence $f\mapsto X_f$ is a Lie algebra homomorphism. One has the identity:
\begin{equation} \label{Jac}
[f,g]=d\lambda(X_f,X_g)+f\, dg(v)-g\, df(v).
\end{equation}
One also  has a projection $TM\to \xi$ along the direction of the Reeb field $v$. Denote by ${\hat u}$ the ``horizontal" part of $u\in TM$, that is, its projection to the contact hyperplane. Then one has: $X_f=fv+{\hat X_f}$.

As in Lemma \ref{Reeb}, assume that the Reeb field $v$ is tangent to a co-Legendrian foliation ${\cal F}$ on a contact manifold $M^{2n-1}$ with a contact form $\lambda$. Let $f_i:M\to \R,\ i=1,\dots,n-1$, be locally defined functions whose common level surfaces are the leaves of ${\cal F}$, and let $u_i=X_{f_i},\ i=1,\dots, n-1$.

\begin{lemma} \label{commute}
The vector fields $v,u_1,\dots, u_{n-1}$ pairwise commute and span the foliation ${\cal F}$.
\end{lemma}

\proof
Let $f$ be an ``integral" of the foliation ${\cal F}$, that is, a function constant on the leaves. The formulas $\lambda(X_f)=f$ and $L_{X_f}\lambda=df(v) \lambda$, along with the Cartan formula,  imply that 
$i_{X_f} d\lambda =df(v) \lambda -df$. It follows that, for every test vector $w\in T{\cal G}$, one has $ d\lambda(X_f,w)=0$: indeed, $\lambda(w)=0$ since $w\in \xi$, and $df(w)=0$ since $w\in T{\cal F}$. It follows that ${\hat X_f}$ lies in the symplectic orthogonal complement to $T{\cal G}$ in $\xi$. Since ${\cal G}$ is Legendrian, ${\hat X_f}$ is tangent to ${\cal G}$. Therefore $X_f$ is tangent to ${\cal F}$.

Next, we claim that $[f_i,f_j]=0$. Indeed, since $v$ is tangent to ${\cal F}$, one has $df_i(v)=0$ for all $i$. It follows from (\ref{Jac}) that $[f_i,f_j]=d\lambda(u_i,u_j)$. Since $v\in\ker d\lambda$, the latter is equal to $d\lambda({\hat u_i},{\hat u_j})$, and this is zero since all ${\hat u_i}$ lie in the Legendrian space $T{\cal G}\subset \xi$.

Likewise, $[1,f_i]=d\lambda(v,u_i)=0$ since $v\in\ker d\lambda$. It follows that the vector fields $v,u_1,\dots, u_{n-1}$ pairwise commute, as claimed.
\proofend

Suppose that a system of differential equations is given. To solve the system in quadratures means to obtain its solution by a finite number of ``algebraic" operations (including inversion of functions) and ``quadratures", integration of known functions, see, e.g., \cite{Kozlov}. 

\begin{theorem} \label{quadratures}
A continuous-time completely integrable contact system is solvable in quadratures.
\end{theorem}

\proof If the contact structure and $n-1$ first integrals $\{ f_i\}$ are given then one can find the special contact as  $\lambda=\lambda_0/\lambda_0(w)$ where $\lambda_0$ is a contact form  and $w$ is a  contact vector field. Then one can compute the commuting vector fields $v,u_1,\dots, u_{n-1}$ from Lemma \ref{commute}. 
The contact vector field defining the dynamical system is a linear combination of these commuting fields with constant coefficients.
 It remains to refer  to a theorem of S. Lie that if $X_1,\dots,X_n$ are commuting and linearly independent vector fields in a domain in $\R^n$ then the differential equation $\dot x=X_1(x)$ is solvable in quadratures, see \cite{Kozlov}.
\proofend

\begin{remark}
{\rm The above theorem is a manifestation of a general phenomenon that the existence of an (explicit) $\R^n$-action 
(and even the semi-direct product action) on a manifold implies  solvability in quadratures, see \cite{Kozlov, Kozlov2} and  references therein.

Note that  the semi-direct product action, defined via the foliation $({\cal F, G})$, does not depend on a contact form, while the definition of the $\R^n$-action requires the knowledge of the special 1-form.}
\end{remark}

%%%%%%%%%%%%%%%%%%%%%%%%%%%%%%%%%%%%%%%%%%%%%%%%%%%%%%%%%%%%%%

\subsection{Example: integrable geodesic flow on a Riemannian manifold} \label{intgeo}

The following is a familiar example from Riemannian geometry.

Let $N$ be a Riemannian manifold of dimension $n$, $T^*N$ its cotangent bundle, $H:T^*N\to\R$ the energy function: $H(q,p)=|p|^2/2$, where $p$ is the momentum and the norm is the Riemannian one. The geodesic flow on $T^*N$ is the Hamiltonian vector field of the function $H$ with respect to the canonical symplectic structure of the cotangent bundle.

Note that $T^*N$ with the zero section deleted is the symplectization of the space $M=ST^*N$ of oriented contact elements in $N$. The homogeneous of degree one Hamiltonian $\sqrt{2H}=|p|$ defines a contact vector field in the contact manifold $M$; this is the geodesic flow on the space of contact elements. The Riemannian metric provides a section of the symplectization $P=T^*N - N \to ST^*N=M$ and hence a contact form on $M$; namely, $M$ is identified with the hypersurface $H=1$. The geodesic flow on $T^*N$ being restricted to $M$ becomes the Reeb vector field of this contact form. 

Assume that the geodesic flow on $T^*N$ is completely integrable: there exist almost everywhere independent and Poisson commuting homogeneous functions $f_1,\dots,f_{n-1}:T^*N\to\R$, invariant under the flow of sgrad $H$. Restricting to the hypersurface $M=\{H=1\}$, one has a co-Legendrian foliation ${\cal F}$, defined by the integrals $f_i$, and the respective Legendrian foliation ${\cal G}$, spanned by the Hamiltonian vector fields sgrad $f_i$.  The Reeb field is tangent to ${\cal F}$, which takes us to the situation of Section \ref{contf}. Thus this geodesic flow is a completely integrable continuous time contact dynamical system.

 The example of this section can be generalized as follows. Let $M$ be a contact manifold, ${\cal F}$ and ${\cal G}$ a co-Legendrian and the Legendrian foliations. In the notation of Section \ref{interp}, assume that $H:P\to\R$ is a homogeneous function  of degree one (replacing energy by the norm), which Poisson commutes with the functions $\bar f_i$. Then the level hypersurface $\{H=1\}$ is a section of the bundle $\pi :P\to M$, and we identify $M$ with this section. 

Recall that $P$ has the Euler field $E$, the symplectic structure $\omega$ and the 1-form $\lambda$ satisfying relations (\ref{sy}). Let $v={\rm sgrad}\ H$.

\begin{lemma} \label{sect} 
The vector field $v$ is the Reeb field of the form $\lambda$.
\end{lemma}

\proof
One has:
$$
i_v \omega = -dH=0, \quad \lambda(v)=\omega (v,E)=dH(E)=H=1,
$$
the first equality due to the fact that $H$ is 1 on the section, and the second to the fact that $H$ is homogeneous of degree one.
\proofend

%%%%%%%%%%%%%%%%%%%%%%%%%%%%%%%%%%%%%%%%%%%%%%%%%%%%%%%%%%%%%%

\subsection{Previous work and generalizations} \label{earlier}
As we already mentioned, contact complete integrability was studied earlier by a number of authors. Here we briefly survey these works.

P. Liberman \cite{Lib} studied Legendrian foliations of contact manifolds endowed with a contact form $\lambda$.  Such a foliation, ${\cal G}$, is called {\it $\lambda$-complete} if the Jacobi bracket of two integrals of  ${\cal G}$ is again an integral (this does not exclude constants). This assumption implies that there exists a flag of foliations $({\cal F}, {\cal G})$ where ${\cal F}$ is co-Legendrian and  tangent to the Reeb field of the form $\lambda$. It is also proved in \cite{Lib} that, in this case, the leaves of ${\cal F}$ and the leaves of ${\cal G}$ have affine structures. In our terms, the former is a consequence of the local $\R^n$-action  by contactomorphisms described in Lemma \ref{commute}, and the latter is a particular case of Lemma \ref{affstr}. Independently, Pang obtained similar results in \cite{Pang}.

In terms of the present paper, A. Banyaga and P. Molino \cite{B-M1,B-M2,Ba} define a  completely integrable contact manifold as a co-Legendrian foliation whose leaves are the orbits of an abelian Lie algebra $\mathfrak{g}$ of contact vector fields. (This point of view was also taken in \cite{Mir}.)
It is proved in \cite{B-M2} that there exists a $\mathfrak{g}$-invariant contact form (note that no assumption on compactness of the respective group of contactomorphisms is made), and that the Reeb field of this contact form belongs to $\mathfrak{g}$. Thus one has the situation of Section \ref{contf}. 

E. Lerman \cite{Ler} studied contact toric manifolds, that is, contact manifolds $M^{2n-1}$ with an action of a torus $\T^n$ by contactomorphisms. This is analogous to the much better studied theory of symplectic toric manifolds, see, e.g., \cite{Aud}.
Let us emphasize that the papers \cite{Lib,Pang,B-M1,B-M2,Ba, Ler, Mir} contain many other interesting  results; we have mentioned only what is relevant to the present work.

\begin{remark} 
{\rm
We also note that weakly integrable contact systems discussed in Section \ref{weaksect}
allow various generalizations to manifolds with more general non-integrable distributions.
Consider a manifold $M$ with a non-integrable distribution $\tau$, which is not necessarily contact. 
One way to define a dynamical system is to consider on such a manifold the action
of a semi-direct product group $\R^l\ltimes \R^k$, which will lead to solvability in quadratures.
Another way is to consider a foliation trasversal to the distribution $\tau$, whose intersections with this distribution have a natural $\R^k$-action. This is the case, in several examples of non-holonomic mechanics, including the non-holonomic oscillator and the Chaplygin skate, which exhibit a weak form of integrability, see \cite{B-C}.
}
\end{remark}

%%%%%%%%%%%%%%%%%%%%%%%%%%%%%%%%%%%%%%%%%%%%%%%%%%%%%%%%%%%%%%
%%%%%%%%%%%%%%%%%%%%%%%%%%%%%%%%%%%%%%%%%%%%%%%%%%%%%%%%%%%%%%

\section{Null lines and the billiard ball map} \label{null}

\subsection{Contact space of oriented light-like lines} \label{paradise}

The space of oriented lines $M^{2n}$ in $\R^{n+1}$  has a canonical symplectic structure, which can be defined as follows  (see, e.g., \cite{A-G}). Start with the cotangent bundle $T^*\R^{n+1}$, and consider the unit energy hypersurvace $|p|^2=2$. The restriction of the canonical symplectic structure on $T^*\R^{n+1}$ to this hypersurface has a one-dimensional kernel. The integral curves of this field of kernels are called the characteristics. A characteristic consists of unit covectors  whose foot points belongs to a fixed line and whose kernels are orthogonal to this line and agree with its orientation. The space of characteristics is again symplectic, and it is identified with the space ${M}$ of oriented lines in $\R^{n+1}$. 

This construction is called symplectic reduction. Symplectic reduction also applies to the space of oriented non-parameterized geodesics of a Riemannian or Finsler manifold (assuming this space is a smooth manifold, which is always  the case locally).

Consider now pseudo-Euclidean space $\R^{p,q}$ with $p+q=n+1$. There are three types of lines: space-like, time-like, and light-like, depending on whether the energy  $|p|^2/2$ is positive, negative or null. Denote these spaces by $M_+^{2n}, M_-^{2n}$ and $M_0^{2n-1}$, respectively. Symplectic reduction on the energy levels $\pm 1$ yields symplectic structures on spaces $M_{\pm}$, but the symplectic reduction on zero energy level yields a space $P^{2n}$ which is different from $M_0$: the condition $|p|^2=0$ still  allows to multiply $p$ by a non-zero real.   $P$ is the space of {\it scaled} null geodesics which fibers over $M_0$ with fiber $\R^*$. Thus $M_0$ is a contact manifold whose symplectization is $P$, the symplectic reduction of $T^*\R^{p,q}$ on zero energy level, see \cite{Kh-Tab} for details.

The space of oriented light-like geodesics was studied about 30 years ago by Yu. Manin in his work on application of twistors to the Yang--Mills equation. Manin called this space {\it paradise} (because it consists of  {\it celestial spheres}, the world lines of photons emanating from  point sources in the Minkowski space $\R^{1,3}$).

%%%%%%%%%%%%%%%%%%%%%%%%%%%%%%%%%%%%%%%%%%%%%%%%%%%%%%%

\subsection{Billiard ball map and accelerating orbits} \label{billball}

The billiard dynamical system in a Riemannian manifold with a smooth boundary describes the motion of a free mass-point (``billiard ball"). The point moves along a geodesic with constant energy until it hits the boundary where the elastic reflection occurs: the normal component of the velocity instantaneously changes sign whereas the tangential component remains the same. This is the billiard flow, a continuous-time system.  The {\it billiard ball map} $T$ acts on oriented geodesics and takes the incoming trajectory of the  billiard ball to the outgoing one. $T$ preserves the symplectic  structure on the space of oriented geodesics. We refer to \cite{Ta1,Ta2} for  information about billiards.

This description applies equally well to billiards in pseudo-Riemannian manifolds, in particular, pseudo-Euclidean spaces. 
A new feature is that now the normal vector to the boundary of the billiard table may be tangent to the boundary; the billiard reflection is not defined at such points. $T$ preserves the type of a billiard trajectory, space-, time-, or light-like. On the spaces $M_{\pm}$, the billiard ball map is still symplectic, but on the space $M_0$, it is a contact transformation, see \cite{Kh-Tab}. 

In fact, we also have a billiard transformation ${\overline T}:P\to P$ on the space of scaled light-like lines described by the reflection law in the opening paragraph of this section. For the projection $\pi:P\to M_0$, one has a commutative diagram: $\pi \circ {\overline T}=T\circ \pi$. 

\begin{figure}[hbtp]
\centering
\includegraphics[width=1.8in]{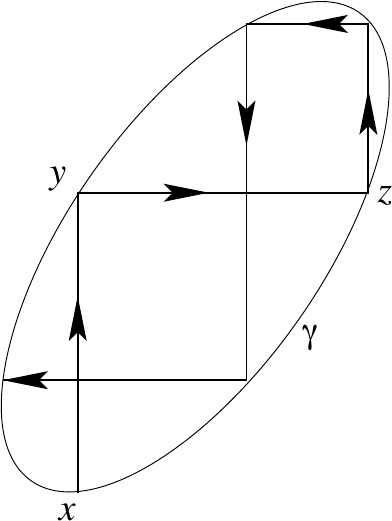}
\caption{A map of an oval}
\label{oval}
\end{figure}

\begin{example} \label{cycle}
{\rm The simplest example is the billiard inside a convex smooth closed curve $\gamma$ in the Lorentz plane $\R^{1,1}$. There are two null directions, say, horizontal and vertical, and the billiard system, restricted to the null directions, is the following self-map of $\gamma$: choose a point $x\in \gamma$, draw the vertical line through $x$ until its second intersections with $\gamma$ at  point $y$, draw the horizontal line through $y$ until its second intersection with $\gamma$ at point $z$, etc., see figure \ref{oval}. This map was studied in various contexts, see \cite{Gen-Kh-Tab} for references.
}
\end{example}

Let us describe an interesting feature of this billiard system, absent in the Euclidean case. Suppose we have a closed  light-like billiard trajectory. Is it possible that, traversing this trajectory, the billiard ball returned to the original position with a different velocity vector, say, greater than the original one? Let us call such a conjectural light-like periodic orbit {\it accelerating}. 

Assume that the metric is $dxdy$, so the null directions are vertical and horizontal. Let $P_1,\dots,P_{2n}\in \gamma$ be the consecutive reflection points of a periodic light-like billiard trajectory, and let $t_i$ be the slope of the curve $\gamma$ at point $P_i$. Consider the billiard ball starting at $P_1$ with, say, unit horizontal velocity, $(1,0)$. Then it will return to point $P_1$ with velocity $(v,0)$.

\begin{lemma} \label{accel}
One has: 
$$v=\frac{t_2 t_4\dots t_{2n}}{t_1 t_3\dots t_{2n-1}}.$$
Further, $v=1$ if and only if the periodic light-like trajectory is stable in the linear approximation.
\end{lemma}

\begin{figure}[hbtp]
\centering
\includegraphics[width=2in]{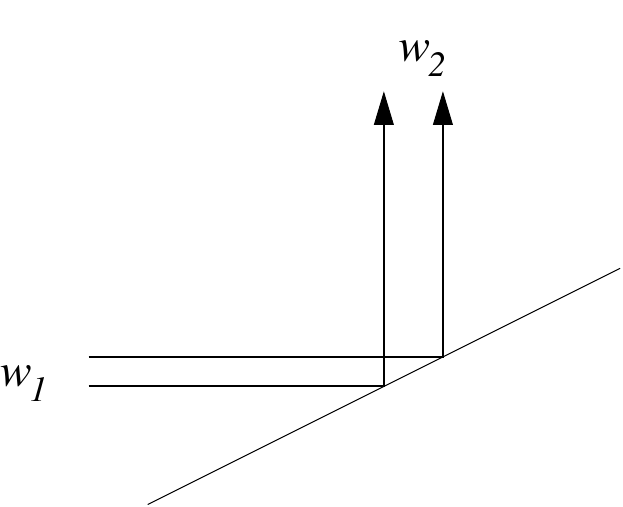}
\caption{Reflection at a point}
\label{refl}
\end{figure}

\proof
Consider an instance of reflection, see figure \ref{refl}. If the slope of $\gamma$ at the reflection point is $t$ then the 
tangent vector to $\gamma$ is $(1,t)$, and the normal vector is $(1,-t)$. Then the reflection is as follows:
$$
(1,0)=\frac{1}{2}(1,t)+\frac{1}{2}(1,-t)\mapsto \frac{1}{2}(1,t)-\frac{1}{2}(1,-t)=(0,t).
$$
Likewise, the vertical-to-horizontal reflection scales the speed down by the slope. This implies the first claim of the lemma.

For the second claim, consider an infinitesimal horizontal beam reflecting in $\gamma$, see figure \ref{refl} again. Let $w_1$ and $w_2$ be the widths of the incoming and the outgoing beams. Then, from elementary geometry, $w_1/w_2=t$. It follows that the condition for the width of the beam to remain the same after all $2n$ reflections is $v=1$. But the former is the linear stability condition for the periodic orbit. 
\proofend

Since the slopes of $\gamma$ at points $P_i$ can be deformed at will (which does not affect the reflection, since the orbit is formed by the same null segments), one can easily construct a billiard table with an accelerating light-like periodic orbit. For such a billiard, there exists no section of the symplectization $\pi:P\to M_0$, invariant under the map ${\overline T}: P\to P$.

\begin{remark} \label{geotoo}
{\rm A similar acceleration phenomenon is possible for closed light-like geodesics on pseudo-Riemannian manifolds.
}
\end{remark}

%%%%%%%%%%%%%%%%%%%%%%%%%%%%%%%%%%%%%%%%%%%%%%%%%%%%%

\subsection{Billiard inside an ellipsoid} \label{billell}

In this section, we turn to our main example, the billiard system inside an ellipsoid in pseudo-Euclidean space. This system was studied in detail in \cite{Kh-Tab}; below we summarize the relevant results.

\begin{figure}[hbtp]
\centering
\includegraphics[width=3in]{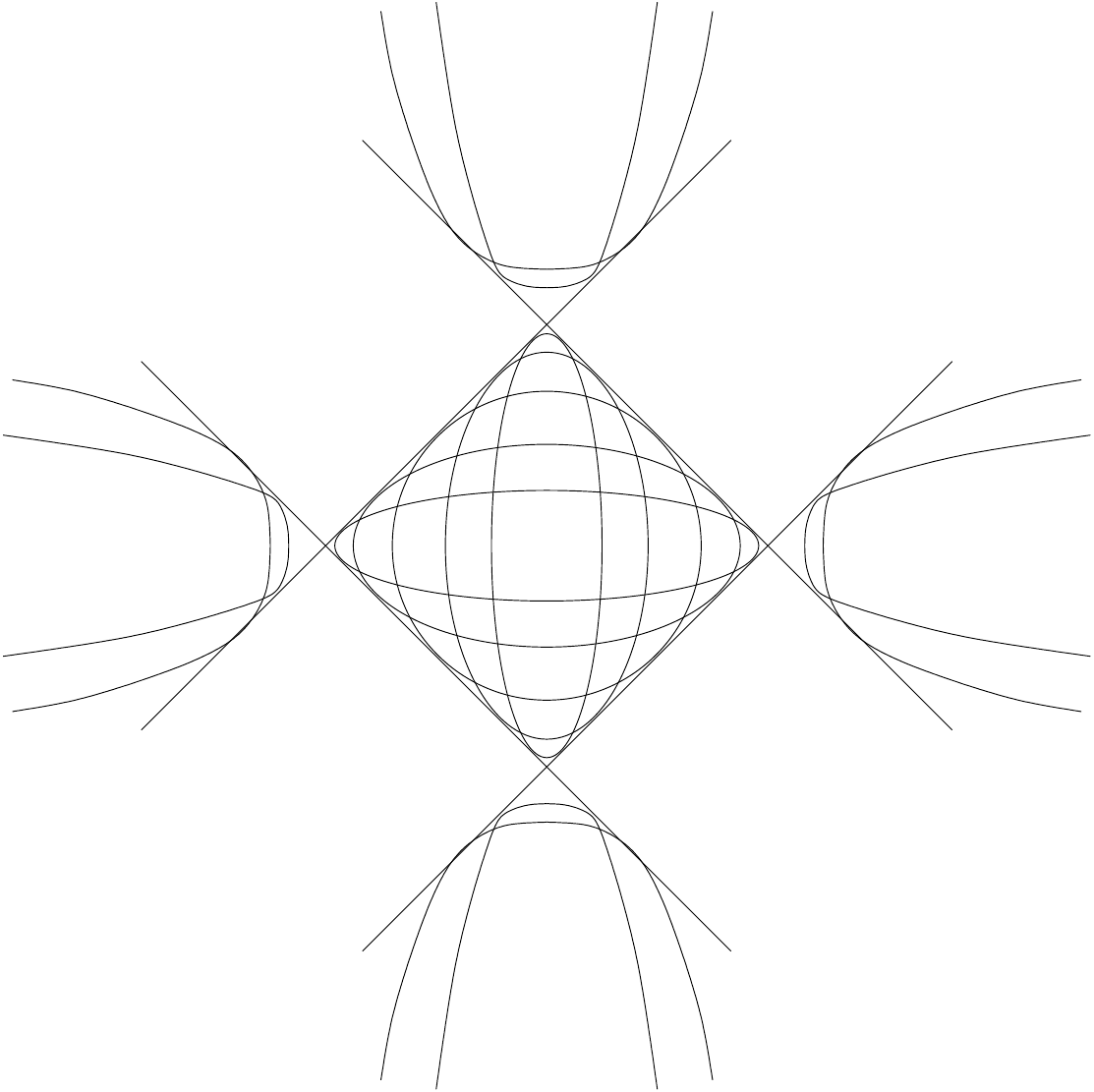}
\caption{A family of pseudo-confocal conics; null directions have slopes $\pm 1$}
\label{domains}
\end{figure}

Consider pseudo-Euclidean space $V^{n+1}=\R^{p,q}$ with $p+q=n+1$, and let $E: V\to V^*$ be the self-adjoint operator such that the metric is given by $E(x) \cdot x$ where dot denotes the pairing between vectors and covectors. Let $A: V\to V^*$ be a positive-definite self-adjoint operator defining an ellipsoid $A(x)\cdot x=1$. Since $A$ is positive-definite, both forms can be  simultaneously reduced to principle axes, and we
assume that $A={\rm diag}(a_1^{-2},\dots,a_n^{-2})$ and $E={\rm diag}(1,\dots,1,-1,\dots,-1)$. Consider the {\it pseudo-confocal family of quadrics} 
\begin{equation} \label{psconf}
\frac{x_1^2}{a_1^2+\lambda} + \frac{x_2^2}{a_2^2+\lambda} +\dots + \frac{x_p^2}{a_p^2+\lambda} + \frac{x_{p+1}^2}{a_{p+1}^2-\lambda}+\dots + \frac{x_{p+q}^2}{a_{p+q}^2-\lambda}=1
\end{equation}
where $\lambda$ is a real parameter (see figure \ref{domains} for a two-dimensional example).
Let $M^{2n-1}_0$ be the contact space  of oriented null lines in $V^{n+1}$, and let $P^{2n}$ be its symplectization, the space of scaled null lines.

The following theorem is proved in \cite{Kh-Tab}.

\begin{theorem} \label{JCh}
1) The tangent lines to a fixed  light-like geodesic on an ellipsoid in  pseudo-Euclidean space $V^{n+1}$ are tangent to  $n-2$ other fixed quadrics from the pseudo-confocal family (\ref{psconf}).\\
2) A  light-like billiard trajectory inside an ellipsoid in pseudo-Euclidean space  $V^{n+1}$ remains tangent to $n-1$ fixed  pseudo-confocal quadrics.\\
3)  The set $N$ of oriented light-like  lines, tangent to   fixed $n-1$  pseudo-confocal quadrics, is a codimension $n-1$  submanifold in  $M_0$, foliated by  Legendrian in $M_0$ submanifolds, which are of codimension one in $N$. 
\end{theorem}

\noindent (For space- and time-like lines, the number of pseudo-confocal quadrics in  statements 1 and 2 is one greater.)

In  terms of the present paper, the null lines, tangent to  $n-1$  fixed  pseudo-confocal quadrics, constitute the leaves of a co-Legendrian foliation ${\cal F}$ in $M_0$. Let $f_1,\dots, f_{n-1}: M_0\to\R$ be smooth functions defining the foliation ${\cal F}$ (these functions index the pseudo-confocal quadrics tangent to a given line), and let $\bar f_i$ be their lifts to $P$, the space of scaled light-like lines.  Then the functions $\bar f_i$ Poisson commute.  

Let us also describe the leaves of the Legendrian foliation ${\cal G}$. A leaf of ${\cal F}$ consists of null lines $\ell$ tangent to fixed $n-1$ pseudo-confocal quadrics, say, $Q_1,\dots, Q_{n-1}$. Let $v_i$ be the geodesic vector field on $TQ_i$. Considering the oriented tangent lines to a geodesic curve, we view $v_i$ as a vector field on the space of lines tangent to $Q_i$. Then these vector fields commute, and the leaf of the foliation ${\cal G}$ through point $\ell$ is generated by the fields $v_1,\dots, v_{n-1}$.

Explicit formulas for integrals are as follows, cf. \cite{Mo1,Mo2} in the Euclidean case.
Identify the tangent $TV$ and cotangent $T^*V$ spaces via the pseudo-Euclidean metric. Then one has the following integrals of the billiard flow on $TV$:
\begin{equation} \label{Moser}
F_k=\frac{v_k^2}{e_k} +\sum_{i\neq k} \frac{(x_iv_k-x_kv_i)^2}{e_i a_k^2-e_k a_i^2},\ \ k=1,\dots,n+1,
\end{equation}
where $x_i$ are the coordinates of the position and $v_i$ of the velocity vectors, and where $e_1=\dots=e_p=1,\ e_{p+1}=\dots=e_{p+q}=-1$. These integrals Poisson commute and satisfy the relation $\sum F_k = \langle v,v\rangle$.  The same formulas give integrals of the geodesic flow on a quadric in pseudo-Euclidean space $V$.
Note that the integrals (\ref{Moser}) are quadratic in velocities.

In the Euclidean case, when all $e_i=1$, the functions $F_k/\langle v,v\rangle$ descend to the space of oriented lines and are integrals of the billiard ball map. In the pseudo-Euclidean case, $\langle v,v\rangle=0$ for the null directions, and one cannot divide by $\langle v,v\rangle$.

Following \cite{Ta2}, let us describe another integral of the billiard ball map, homogeneous of degree one in the velocity. Let $x$ be a point of the ellipsoid and $v$ an inward vector with foot point $x$. As before, one has the billiard ball  transformation ${\overline T}$ on such tangent vectors. If $v$ is null then the set of the inward tangent vectors with foot point on the ellipsoid is identified with the space of scaled oriented lines $P$. 

\begin{proposition} \label{newint}
1) The function $H(x,v):=Ax \cdot v$ is negative.\\
2) $H(x,v)$ is invariant under the billiard ball  transformation ${\overline T}$.\\
3) $H(x,v)$ Poisson commutes  with the functions $\bar f_i,\ i=1,\dots, n-1$. 
\end{proposition}

\proof For the first claim, note that  $Ax$ is the outward normal covector and $v$ has the inward direction, hence $H(x,v)<0$.

\begin{figure}[hbtp]
\centering
\includegraphics[width=3in]{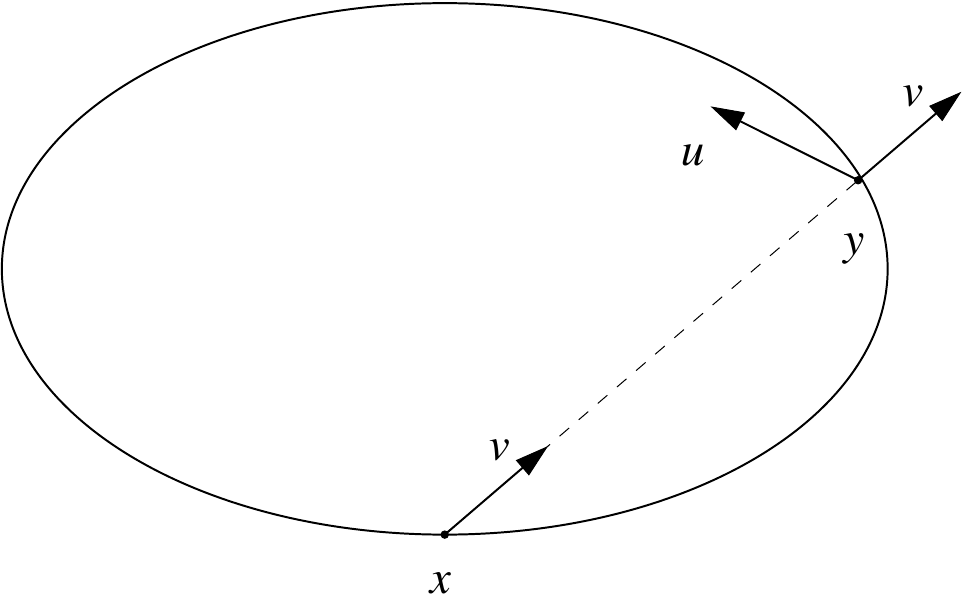}
\caption{Billiard ball map as the composition of two involutions}
\label{reflection}
\end{figure}

For the second claim, the billiard ball map is the composition of two maps: $(x,v)\mapsto (y,v) \mapsto (y,u)$, where the second is the billiard reflection, see figure \ref{reflection}. 
We claim that  $Ax \cdot v =- Ay \cdot v = Ay \cdot u$.

To prove the first equality, note that $(Ax + Ay) \cdot (y-x)=0$ since $A^*=A$ and $Ax\cdot x= Ay\cdot y=1$. On the other hand, $v$ is collinear with $y-x$, hence $Ax \cdot v =- Ay \cdot v$.

To prove the second equality, note that, due to the reflection law, $v+u$ is a tangent vector to the ellipsoid at point $y$. On the other hand, $Ay$ is the normal covector to the ellipsoid $Ay \cdot y=1$. Hence $Ay \cdot v = -Ay \cdot u$, as claimed.

Now we prove the third claim. Extend the function $H$ to the tangent bundle $TV$ in such a way that it is invariant along straight lines: $H(x+tv,v)=H(x,v)$ for all $t\in\R$. Then this extended function is an integral of the billiard flow inside the ellipsoid since it is also invariant under the reflection, see above. Since $F_k$ is a complete system of first integrals, $H$  is functionally dependent on  integrals $F_k$ in (\ref{Moser}). Note that each $F_k$ is also invariant  along straight lines: $F_k(x+tv,v)=F_k(x,v)$. Hence the functional relation descends, in particular, to the space of scaled null lines $P$.

Therefore, it suffices to show that the functions $\bar f_i$ and $F_k$ Poisson commute in the space $P$ of scaled light-like lines. 
Indeed, as we mentioned above, the Hamiltonian vector field $v_i = {\rm sgrad} \bar f_i$ defines the geodesic flow on the pseudo-confocal quadric $Q_i$, and the functions $F_k$ are integrals of the geodesic flow on these quadrics, hence $\{\bar f_i, F_k\}=0$.
\proofend

Proposition \ref{newint} places us in the situation of Section \ref{intgeo}, and therefore, of Section \ref{commf}: one has an invariant contact form on the space of null geodesics whose Reeb field is tangent to the co-Legendrian foliation. More specifically, the integral $H(x,v)$ provides a section of the symplectization bundle $P\to M_0$, and hence a special contact form on $M_0$. It follows that the billiard ball map on the light-like oriented lines inside an ellipsoid in pseudo-Euclidean space is a completely integrable contact transformation.

\begin{remark} \label{lim}
{\rm One can view the contact integrable system on null geodesics on an ellipsoid as a limit of the Hamiltonian completely integrable system of space-like geodesics. While only $n-1$ independent integrals of the Hamiltonian system (out of the $n$) survive in the limit, when passing to the contact manifold, the corresponding $\R^n$-action on the symplectic manifold of space-like geodesics does extend to the space of null geodesics on the ellipsoid. 
}
\end{remark}

%%%%%%%%%%%%%%%%%%%%%%%%%%%%%%%%%%%%%%%%%%%%%%%%%%%%
%%%%%%%%%%%%%%%%%%%%%%%%%%%%%%%%%%%%%%%%%%%%%%%%%%%%

\bigskip

{\bf Acknowledgments}. We are grateful to M. Audin, A. Banyaga, E. Lerman, Yu. Manin, I. Miklashevskii and V. Ovsienko  for stimulating discussions. B.K. is grateful to the MSRI for kind hospitality during the work on this paper.  The  authors were partially supported by NSERC and NSF grants, respectively.


\bigskip



\begin{thebibliography}{99}

\bibitem{Arnold}
V. Arnold. {\it Arnold's problems}. Springer-Verlag, Berlin, Heidelberg, and Phasis, Moscow, 2004. 

\bibitem{A-G} V. Arnold, A. Givental. {\it Symplectic geometry}, 1-136. Encycl. of Math. Sci., Dynamical Systems, 4, Springer-Verlag, 1990.

\bibitem{Kozlov} V.I.Arnold, V.V.Kozlov, A.I.Neishtadt. {\it Mathematical aspects of classical and celestial mechanics}, Encycl. of Math. Sci., Dynamical systems 3. Springer-Verlag, 2006. xiv+518 pp.

\bibitem{Aud} M. Audin. {\it Torus actions on symplectic manifolds.} Birkhauser Verlag, Basel, 2004.

\bibitem{B-M1} A. Banyaga, P. Molino. {\it G\'eom\'etrie des formes de contact compl\'etement int\'egrables de type toriques.} S\'eminaire Gaston Darboux de G\'eom\'etrie et Topologie Diff\'erentielle, 1991--1992, 1--25, Univ. Montpellier II, Montpellier, 1993.

\bibitem{B-M2} A. Banyaga, P. Molino. {\it Complete integrability in contact geometry.} Penn State preprint PM 197, 1996.

\bibitem{Ba} A. Banyaga. {\it The geometry surrounding the Arnold-Liouville theorem.} Advances in geometry, 53--69, Progr. Math., 172, BirkhŠuser Boston, Boston, MA, 1999.

\bibitem{B-C} L. Bates, R. Cushman. {\it What is a completely integrable nonholonomic dynamical system?}
Reports on Math. Physics, {\bf 44} (1999), no.1/2, 29--35.

%\bibitem[Lichnerowicz and K]{DLM} P. Dazord, A. Lichnerowicz, Ch.-M. Marle. \textit{Structure locale des vari\'et\'es de Jacobi.} J. Math. Pures Appl. (9) {\bf 70} (1991),  101--152. 

\bibitem{Gen-Kh-Tab} D. Genin, B. Khesin, S. Tabachnikov.
{\it Geodesics on an ellipsoid in Minkowski space.} Enseign. Math.  {\bf 53} (2007),  307--331.

\bibitem{Kh-Tab} B. Khesin, S. Tabachnikov. {\it Pseudo-Riemannian geodesics and billiards.} Adv. in Math., {\bf 221} (2009), 1364-1396.

\bibitem{Kob} S. Kobayashi. {\it Transformation groups in differential geometry.}  Springer-Verlag, Berlin, 1995.

\bibitem{Kozlov2} V. Kozlov. {\it Integrability and non-integrability in Hamiltonian mechanics.}  Russ. Math. Surv. {\bf 38:1} (1983) 1--76.

\bibitem{Ler} E. Lerman. {\it Contact toric manifolds.}  J. Symplectic Geom.  {\bf 1}  (2003),  785--828.

\bibitem{L-T} M. Levi, S. Tabachnikov. {\it The Poncelet grid and billiards in ellipses. }
Amer. Math. Monthly {\bf 114} (2007),  895--908. 

\bibitem{Lib} P. Libermann. {\it Legendre foliations on contact manifolds.} 
Differential Geom. Appl. {\bf 1} (1991),  57--76. 

\bibitem{Mir}   E. Miranda. {\it A normal form theorem for integrable systems on contact manifolds.}  Proceedings of the XIII Fall Workshop on Geometry and Physics, Publ. R. Soc. Mat. Esp., {\bf 9}  (2005), 240--246.

%{\it On symplectic linearization of singular Lagrangian foliations.}  PhD thesis, June, 2003, http://www.picard.ups-tlse.fr/~miranda/research.html

%\bibitem[Monnier-Zung]{MZ} Ph. Monnier, N.T. Zung. {\it Normal forms of vector fields on Poisson manifolds}.  Ann. Math. Blaise Pascal  {\bf 13}  (2006),  no. 2, 349--380.

\bibitem{Mo1}  J. Moser. {\it Various aspects of integrable Hamiltonian systems}, 233-289. Progr. in Math. {\bf 8}, Birhauser, 1980. 

\bibitem{Mo2} 
J. Moser. {\it Geometry of quadrics and spectral theory}, 147-188.  Chern Symp. 1979,
Springer-Verlag, 1980.

\bibitem{Pang} M-Y. Pang. {\it The structure of Legendre foliations. }
Trans. Amer. Math. Soc. {\bf 320} (1990),  417--455. 

\bibitem{Ta1} S. Tabachnikov. {\it Billiards}. Soc. Math. de France,
Paris, 1995.

\bibitem{Ta2} S. Tabachnikov. {\it Geometry and billiards}.  Amer. Math. Soc., Providence, RI, 2005.

%\bibitem[Toth-Zelditch]{Toth-Zel} J.A. Toth, S. Zelditch, {\it Riemannian manifolds with uniformly bounded eigenfunctions.}  Duke Math. J.  111  (2002),  no. 1, 97--132, arXiv.org:math-ph/0002038

\bibitem{Ve1}  A. Veselov. {\it Integrable mappings.}  Russian Math. Surveys {\bf 46} (1991), no. {5}, 1--51.
 

\end{thebibliography}
\end{document}